\documentclass[10pt]{article}
\usepackage{graphicx}

\usepackage{xcolor,paralist,hyperref,titlesec,fancyhdr,etoolbox}
\usepackage{amsmath,amsfonts,amssymb,amsthm,graphicx}
\usepackage[mathscr]{eucal}
\usepackage{enumerate}
\usepackage{mathrsfs}
\usepackage{multicol}
\newtheorem{defn}[subsection]{Definition}
\newtheorem{thm}[subsection]{Theorem}
\newtheorem{cor}[subsection]{Corollary}

\newtheorem{prop}[subsection]{Proposition}

\newtheorem{eg}[subsection]{Example}

\usepackage{lipsum}
\usepackage{authblk}
\usepackage{blindtext}
\begin{document}
\title{Direct Sum of Lower Semi-Frames in Hilbert Spaces}

\author{Hemalatha M$^1$, P. Sam Johnson$^1$, Harikrishnan P.K.$^{3}$}
\date{}
\maketitle
$^{1}$ Department of Mathematical and Computational Sciences, \\National Institute of Technology Karnataka,  Mangaluru, 575025, Karnataka, India.\\
\indent $^{3}$ Department of Mathematics, Manipal Institute of Technology,\\  Manipal Academy of Higher Education,  Manipal, 576104, Karnataka, India.
$$email:hemalatham.217ma005@nitk.edu.in,sam@nitk.edu.in,pk.harikrishnan@manipal.edu$$


\begin{abstract}
In this paper, structural properties of lower semi-frames in separable Hilbert spaces are explored with a focus on transformations under linear operators (may be unbounded). Also, the direct sum of lower semi-frames, providing necessary and sufficient conditions for the preservation of lower semi-frame structure, is examined.
\end{abstract} 
AMS Classification Code:42C15, 47A05, 46B15, 40A05.
\section{Introduction}

In a separable Hilbert space, there are several essential classes of sequences that play a pivotal role in the analysis, representation, and reconstruction of elements. These include Bessel sequences, frames, lower semi-frames, and Riesz-Fischer sequences, all of which are interconnected and widely applicable across various disciplines, such as signal processing, functional analysis, and quantum mechanics. The concept of frames was first introduced by Duffin and Schaeffer in 1952 in their study of non-harmonic Fourier series, where they aimed to generalize orthonormal bases. Frames gained significant attention in the early 1990s, Daubechies et al. reintroduced frames in their study on wavelet theory.  With frames, each element of a Hilbert space can still be expressed as a linear combination of the frame elements, even though the elements are not linearly independent. This makes frames especially practical for real-world problems where redundancy, robustness, and flexibility are required. Over time, the concept of frames has been further generalized into different forms to meet the demands of various applications. Notable generalizations include fusion frames, generalized frames, K-frames,
continuous frames, and continuous fusion frames. These generalizations have led to significant advancements in areas such as distributed signal processing, image compression, and quantum information theory. For a detailed exploration of frame theory and its applications, see references \cite{7,K-frame,8}.

In 2011, Casazza, Christensen, and other mathematicians took frame theory a step further by characterizing sequences using their associated frame-related operators. These operators provide a systematic way to understand the structure and properties of frames and related sequences.
An important question in frame theory concerns methods for generating new frames from existing ones. For example, given two frames in a Hilbert space $H$, a natural idea is to consider their sum. However, the sum of two frames does not necessarily result in another frame for $H$. Researchers have addressed this problem by providing conditions under which the sum of two frames forms a frame. Obeidat et al. (2009) and Najati et al. (2013) \cite{5,3} established necessary and sufficient conditions for the sum of two frames to form a new frame, and Young et al.2013, \cite{Dsum} studied about direct sum of frames. Additional results on this topic were explored in \cite{6,5}. Panayappan \cite{tensor}  studied the tensor sum and tensor products of frames in 2014.

In 2012, Antoine and Balazs \cite{1,2} introduced a broader class of sequences called semi-frames. Semi-frames relax the requirement that frames must satisfy both upper and lower bounds. Instead, they require only one bound, making them a more flexible tool for analysis. These semi-frames are also give reconstruction formulae on separable Hilbert spaces. One significant implication of semi-frames is that they allow frame-related operators to be unbounded, which opens up new possibilities for studying frame theory in the context of unbounded operator theory. The motivation of this work is to generate more frame-like sequences using the operator, sum, and direct sum.\\
\indent This paper is arranged as follows. In Section 2, we revise the necessary definitions and results in the literature, which are required in the sequel. Section 3 contains results based on operators that preserve the frame-like properties of the sequences. In Section 4, we discuss the properties of the sum of sequences.  We provide the necessary and sufficient conditions for the direct sum of sequences in Section 5.  
	
\section{Preliminaries}
	In this section, we start by defining the notations used in this paper, revisiting essential fundamental concepts in frames and a few established results that will facilitate understanding the subsequent sections.
    Throughout the paper, we use the notations such as $\mathbb{N}$ be the set of natural numbers. Let $H, H_1, H_2, H_1',H_2'$ be separable Hilbert spaces. We denote $H_1 \oplus H_2 = \{x\oplus y : x \in H_1 \text{ and } y \in H_2\}$ and define $\langle x_1\oplus y_1 , x_2\oplus y_2\rangle = \langle x_1, y_1 \rangle_{H_1}+\langle x_2, y_2 \rangle_{H_2}$. For two subspaces $K_1$ and $K_2$ of $H$, we let $ K_1+K_2$ denote the (set-theoretic) sum of two subspaces, $K_1\oplus K_2$ denotes the internal direct sum of subspaces, where as $K_1\oplus ^\perp  K_2$ denotes the internal orthogonal direct sum of subspaces. For two sequences $F=\{f_n\}_{n\in \mathbb{N}}$ and $G=\{g_n\}_{n\in \mathbb{N}}$, $F \oplus G = \{f_n \oplus g_n:n\in \mathbb{N}\}$. Let $T$ be a linear operator on $H$. The space of all bounded linear operators from $H_1$ to $H_2$ is denoted by $B(H_1, H_2)$.
$Dom(T)$, $R(T)$, $T^*$ denotes the domain of $T$, the range of $T$ and the adjoint of $T$ respectively."  
\begin{defn} \cite{1,2002,7}
		Let $F = \{f_n\}_{n\in \mathbb{N}}$ be a sequence in $H$. 
  \begin{enumerate}[(i)]
 \item  $F$ is said to be a complete sequence for $H$ if $\overline{span\{f_n\}_{n \in \mathbb{N}}} = H$.
      \item $F$ is said to be a frame for $H$ if there exist constants $0< A \leq B < \infty$ such that
		\begin{equation*}
		A\| f \| ^2 \leq \sum_{n \in \mathbb{N}} | \langle f, f_n \rangle |^2 \leq  B\| f \| ^2, \hspace*{0.5cm}  \text{for all     } f \in H.\end{equation*}
  \item $F$ is said to be a Bessel sequence for $H$ if there exists a constant $B>0$ such that
		\begin{equation*}
	\sum_{n \in \mathbb{N}} | \langle f, f_n \rangle |^2 \leq  B\| f \| ^2, \hspace*{0.5cm}  \text{for all     } f \in H.
 \end{equation*}
 
  \item $F$ is said to be a lower semi-frame for $H$ if there exists a constant $A>0$ such that
		\begin{equation*}
		A\| f \| ^2  \leq \sum_{n \in \mathbb{N}} | \langle f, f_n \rangle |^2, \quad  \text{ for all } f \in H.
		\end{equation*}
  $F$ is said to be a lower semi-frame sequence for $H$ if $\{f_n\}_{n \in \mathbb{N}}$ is a lower semi-frame for $\overline{span\{f_n\}_{n \in \mathbb{N}}}$.
 
 \item $F$ is said to be a Riesz-Fischer sequence for $H$ if there exists a constant $A>0$ such that
	\begin{equation*}
	A \sum_{n \in \mathbb{N}} | c_n| ^2 \leq \| \sum_{n \in \mathbb{N}} c_n f_n \| ^2 
	\end{equation*}
	for all finite scalar sequences $\{c_n\}_{n \in \mathbb{N}} \in \ell_2$.
  
   \end{enumerate}
	\end{defn}

	For any given sequence $F = \{f_n\}_{n \in \mathbb{N}} \subseteq H$, we can associate the following three operators :
	\begin{enumerate}[(i)]
				\item The synthesis operator $D: Dom(D) \subseteq \ell^2 \rightarrow H$ is defined by,
		\begin{eqnarray*}
			Dom(D) &=& \Big\{\{c_n\} \in \ell^2 : \sum_{n=1}^{\infty} c_nf_n \text{ is convergent in }H \Big\}\\
			D\{c_n\} &=& \sum_{n=1}^{\infty} c_nf_n,   \text{ for } \{c_n\} \in Dom(D).
		\end{eqnarray*}
		\item The analysis operator $C$ is defined as the adjoint of the synthesis operator $D$. This means that $C=D^*$ is always a closed operator.
		\item The frame operator $S $ is defined as $S = DC."$
		
	\end{enumerate}
 \begin{defn}\cite{kato}
     Let $A$ be a densely defined closed operator on $H$. The reduced minimum modulus for $A$ is $\gamma(A) = \inf \Big\{ \frac{\|Af\|}{\|f\|} : f \neq0 \in Dom(A) \cap N(A)^\perp\Big\}.$
 \end{defn}
 \begin{prop} \cite{4a} \label{closed}
     Let $A$ be a densely defined closed operator on $H$. Then the following statements are equivalent:
     \begin{enumerate}[(i)]
         \item $\gamma(A) >0$.
         \item $R(A)$ is closed.
         \item $R(A^*)$ is closed.
     \end{enumerate}
 \end{prop}

\begin{prop}\label{1} \cite{2011}
	Let $\{f_n\}_{n \in \mathbb{N}} $ be a sequence  in $H$. Then the following statements are valid.
	\begin{enumerate}[(i)]
     \item $\{f_n\}_{n \in \mathbb{N}} $ is a complete sequence for $H$ if and only if $C$ is injective.
		\item   $\{f_n\}_{n \in \mathbb{N}} $ is a lower semi-frame for $H$ if and only if $C$ is injective and $R(C)$ is closed.
		\item $\{f_n\}_{n \in \mathbb{N}} $ is a Riesz-Fischer sequence for $H$ if and only if $C$ is surjective.
        
	\end{enumerate}
\end{prop}
\begin{defn}\cite{Dsum} \label{def}
    Let $F =\{f_n\}_{n \in \mathbb{N}}$ and $G =\{g_n\}_{n \in \mathbb{N}}$ be frames for $H_1$ and $H_2$ respectively.
    \begin{enumerate}[(i)]
        
        \item Let $F'$ and $G'$ be frames for $H_1'$ and $H_2'$ respectively. Then the pair $(F,G)$ is said to be similar to the pair $(F',G')$ if there are invertible $T \in B(H_1, H_1')$ and $S \in B(H_2, H_2')$ such that $T(F)=F'$ and $S(G)=G'$. The pairs are unitarily equivalent if $T$ and $S$ can be taken to be unitary.
        \item $(F, G)$ is said to be a strong complementary pair if there are Parseval frames $F'$ for $H_1'$ and $G'$ for $H_2'$ such that $F' \oplus G'$ is an orthonormal basis of $H_1' \oplus H_2'$ and $(F, G)$ is similar to $(F', G')$
        \item $(F,G)$ is said to be a complementary pair if $F \oplus G$ is a Riesz basis of $H_1 \oplus H_2$.
        \item $(F, G)$ is said to be a strongly disjoint (orthogonal) pair if there are Parseval frames $F'$ for $H_1'$ and $G'$ for $H_2'$ such that $F' \oplus G'$ is a Parseval frame for $H_1' \oplus H_2'$ and $(F ,G)$ is similar to $(F', G')$.
        \item $(F,G)$ is said to be a disjoint pair if $F \oplus G$ is a frame for $H_1 \oplus H_2$.
        \item $(F,G)$ is said to be a weakly disjoint pair if $F \oplus G$ is complete in $H_1 \oplus H_2$.
    \end{enumerate}
\end{defn}
\begin{thm}\cite{5}
	Let $F=\{f_n\}_{n \in \mathbb{N}}$ and $G=\{g_n\}_{n \in \mathbb{N}}$ be frames for $H_1$ and $H_2$ respectively.
\begin{enumerate}[(i)]
    \item $F$ and $G$ are strongly complementary if and only if $\ell^{2}(\mathbb{N})=R(C_1) \oplus ^\perp R(C_2)$.
    \item $F$ and $G$ are complementary if and only if  $\ell^{2}(\mathbb{N})=R(C_1)\oplus R(C_2)$.
    \item $F$ and $G$ are strongly disjoint (orthogonal) if and only if $R(C_1) \perp R(C_2)$.
    \item $F$ and $G$ are weakly disjoint if and only if $R(C_1) \cap R(C_2)=\{0\}.$
\end{enumerate}
\end{thm}
\section{Generating sequences using operators}
 In this section, we explore the properties of operators that preserve the characteristics of general sequences. 
This investigates the behavior of various structured sequences in a Hilbert space $H$, such as complete sequences, lower semi-frames, frames, and Riesz-Fischer sequences, under the action of densely defined linear operators. A fundamental result establishes that any densely defined closed operator on $H$ can be realized as the analysis operator associated with a suitably chosen sequence. Further, it is shown that for a complete sequence $\{f_n\} \subset H$, the image sequence $\{Lf_n\}$ remains complete if and only if the operator $L$ has a dense range. In the case of lower semi-frames, the stability of the structure under the action of an operator depends crucially on the properties of both the operator and its adjoint. Specifically, if $L$ has a dense range and its adjoint $L^*$ has a closed or surjective range, then the image sequence $\{Lf_n\}$ retains the lower semi-frame property. These conditions are also shown to be sufficient for preserving the Riesz-Fischer property. 

Let $\{f_n\}_{n\in \mathbb{N}}$ be a sequence in $H$ with the analysis operator $C$, the synthesis operator $D$ and the frame operator $S$. Let $L: Dom(L) \subseteq H \rightarrow H$ be an operator with dense domain. For the sequence $\{Lf_n\}_{n \in \mathbb{N}},$ the corresponding analysis, synthesis and frame operators are $CL^*$, $LD$ and $LSL^*$ respectively.

\begin{prop}\label{d}
        Let $L$ be a densely defined closed operator on $H$. Then there exists a sequence $\{f_n\}_{n \in \mathbb{N}}$ in $H$ such that $L$ serves as the analysis operator for $\{f_n\}_{n\in \mathbb{N}}$. 
    \end{prop}
    \begin{proof}
       Since $L$ is densely defined closed operator, $L^*$ is densely defined and $(L^*)^* = L.$
        Let $\{e_n\}_{n\in \mathbb{N}}$ be the standard orthonormal basis for $\ell_2$ and define $f_n = L^*e_n$. The analysis operator for $\{f_n\}_{n\in \mathbb{N}}$ is $\{\langle f,f_n\rangle \}_{n\in \mathbb{N}} = \{\langle f,L^*e_n\rangle \}_{n\in \mathbb{N}}
            = \{\langle Lf,e_n\rangle \}_{n\in \mathbb{N}}
            = Lf$ for all $f \in Dom(L)$. Hence, $L$ is the analysis operator for $\{f_n\}_{n \in \mathbb{N}}$.
        \end{proof}
   
\begin{prop}\label{a1}
    Let $\{f_n\}_{n \in \mathbb{N}}$ be a complete sequence in $H,$ and let $L$ be an operator  with dense domain. Then $R(L)$ is dense in $H$ if and only if $\{Lf_n\}_{n \in \mathbb{N}}$ is a complete sequence in $H$.
\end{prop}
\begin{proof}
    The analysis operator $C$ for $\{f_n\}_{n\in \mathbb{N}}$ is injective.
    Since $H = N(L^*) \oplus \overline{R(L)}$ and $R(L)$ is dense, it follows that $ N(L^*) = \{0\}$. 
    Consequently, $L^*$ is injective, which implies that $CL^*$ is also injective. Therefore, $\{Lf_n\}_{n \in \mathbb{N}}$ is a complete sequence in $H$.
    
    Conversely, suppose that $\{Lf_n\}_{n \in \mathbb{N}}$ is a complete sequence in $H$,
    then $CL^*$ is injective.
    If $f \in N(L^*)$, then $ f \in N(CL^*)$.
    Given that $N(CL^*) = \{0\}$, it follows that
          $ f = 0$. Hence, $L^*$ is injective.
    \end{proof}
\begin{prop}\label{a}
    Let $\{f_n\}_{n \in \mathbb{N}}$ be a lower semi-frame for $H$ and let $L$ be an operator  with dense domain on $H$. If $R(L)$ is dense in $H$ and $R(L^*)$ is closed, then $\{Lf_n\}_{n \in \mathbb{N}}$ is also a lower semi-frame for $H$.
\end{prop}
\begin{proof}
    Our assumption on $\{f_n\}_{n \in \mathbb{N}}$ implies that $C$ is injective, closed operator with closed range. Consequently, $C^{-1}$ is bounded on $R(C)$.
     By the Proposition \ref{a1}, $CL^*$ is injective. 
    Now, let $g_n \in R(CL^*)$ such that $g_n \rightarrow g$ as $n \rightarrow \infty$. Then, there exists $h_n \in Dom(CL^*)$ such that $g_n = CL^* h_n$. Since $R(C)$ is closed, there exists $h \in Dom(C)$ such that $CL^* h_n \rightarrow C h.$ Moreover, $L^* h_n \rightarrow h$ because $C^{-1}$ is bounded. As $R(L^*)$ is also closed, there exists $h' \in Dom(L^*)$ such that $L^* h_n \rightarrow L^*h'$. Therefore, $CL^* h_n \rightarrow CL^*h',$ which implies $R(CL^*)$ is closed. Thus, $ \{Lf_n\}_{n \in \mathbb{N}} $ is a lower semi-frame for $H$.    
\end{proof}
 The converse of the Proposition \ref{a} is not true in general. Consider the following example.
\begin{eg} \label{eg 1}
    Let $H = \ell_2$ and $\{f_n\} = \{ne_n\}$, where $\{e_n\}$ is the standard orthonormal basis for $\ell_2$. It is known that $\{f_n\}$ is a lower semi-frame for $\ell_2$. Define an operator $L : Dom(L) \subseteq \ell_2 \rightarrow \ell_2$ as $L(e_n) = \frac{e_n}{n}$. Then, $L$ is densely defined and self-adjoint ($L^* = L$). The sequence $\{(1,0,0,\ldots),(1,\frac{1}{2},0,\ldots),\ldots\}$ converges to $(1,\frac{1}{2},\frac{1}{3},\ldots,\frac{1}{n},\ldots) \notin R(L)$. Therefore, $R(L)$ is not closed. However, the sequence $\{Lf_n\} = \{e_n\}$ forms a frame for $\ell_2$.
\end{eg}
If we modify the sequence in Example \ref{eg 1} to $\{f_n\} = \{e_n\}$ then $\{Lf_n\} =\{\frac{e_n}{n}\},$ which is not a lower semi-frame for $\ell_2$.
\begin{prop}
   Let $\{f_n\}_{n \in \mathbb{N}}$ be a sequence in $H$ with the analysis operator $C$. Let $L$ be an operator  with dense domain on $H$ such that $R(L)$ is dense in $H$ and $R(L^*) =H$. Then $\{Lf_n\}_{n \in \mathbb{N}}$ is a lower semi-frame for $H$ if and only if $\{f_n\}_{n \in \mathbb{N}}$ is a lower semi-frame for $H$.
\end{prop}
\begin{proof}
    The analysis operator for $\{Lf_n\}_{n \in \mathbb{N}}$ is $CL^*$. Since $R(L)$ is dense in $H$, $L^*$ is injective. By our assumption on the sequence $\{Lf_n\}_{n \in \mathbb{N}}$ gives that $CL^*$ is injective, which implies that $C$ is injective on $H$. Since $R(L^*)=H$, $R(CL^*)=R(C)$ is closed. This shows that $\{f_n\}_{n \in \mathbb{N}}$ is a lower semi-frame for $H$. Conversely, the result follows from the Proposition \ref{a}.    
\end{proof}
 \begin{prop}
           Let $\{f_n\}_{n \in \mathbb{N}}$ be a frame for $H$ and $L$ be an operator  with dense domain on $H$. If $\{Lf_n\}_{n \in \mathbb{N}}$ is a lower semi-frame for $H$, if and only if $R(L)$ is dense and $R(L^*)$ is closed.
      \end{prop}
      \begin{proof}
         Let $C$ be the analysis operator for $\{f_n\}_{n \in \mathbb{N}}$. Then $Dom(C) = H$ and $C$ is injective. Since $\{Lf_n\}_{n \in \mathbb{N}}$ is lower semi-frame, $CL^*$ is injective. Hence, this implies that $L^*$ is injective, $R(L)$ is dense. Let $\{g_n\}_{n \in \mathbb{N}} \in R(L^*)$ such that $g_n \rightarrow g\in H$ as $n \rightarrow \infty.$ Since $Dom(C) = H \text{ and } R(CL^*)$ is closed, $\{g_n\}_{n \in \mathbb{N}} \in Dom(C) \text{ and there exists } h\in Dom(CL^*) \text{ such that } Cg_n \rightarrow CL^*h$. Then $g_n \rightarrow L^*h$ because $C^{-1}$ is bounded. Therefore $R(L^*)$ is closed. The converse follows from the Proposition \ref{a}.    
      \end{proof}
\begin{prop}\label{e}
   Let $\{f_n\}_{n \in \mathbb{N}}$ be a lower semi-frame $H$ and $L$ be an operator  with dense domain on $H$. If $R(L^*)$ is closed, then $\{Lf_n\}_{n \in \mathbb{N}}$ is a lower semi-frame sequence for $H$.
\end{prop}
\begin{proof} Since $\overline{span\{Lf_n\}_{n \in \mathbb{N}}} \subseteq \overline{R(L)},$ $R(L) = \overline{span\{Lf_n\}_{n \in \mathbb{N}}} \oplus \overline{span\{Lf_n\}_{n \in \mathbb{N}}}^\perp$. Let $f \in \overline{span\{Lf_n\}_{n \in \mathbb{N}}}^\perp$, which implies $\{\langle f, Lf_n \rangle\}_{n\in \mathbb{N}} = 0$. Then $L^*f = 0$, because $\{f_n\}_{n \in \mathbb{N}}$ is a lower semi-frame. Injectivity of $L^*$ on $\overline{R(L)}$ gives that $f = 0.$ Therefore, $\overline{span\{Lf_n\}_{n \in \mathbb{N}}} = \overline{R(L)}.$ Now consider $\{Lf_n\}_{n \in \mathbb{N}}$ on $\overline{R(L)}$. Then the analysis operator $C_1$ is, $$C_1f = \{\langle f, Lf_n\rangle\}_{n\in \mathbb{N}} = \{\langle L^*f, f_n\rangle\}_{n\in \mathbb{N}} = CL^* \text{ for every }f \in Dom(L^*) \cap \overline{R(L)}.$$  Therefore $C_1 = CL^*|_{Dom(L^*) \cap \overline{R(L)}} $. \\ \hspace*{0.5cm}Assumption on $\{f_n\}_{n\in \mathbb{N}}$ gives that,  $C$ is an injective closed operator with a closed range. Then $C^{-1}$ is bounded. Since $H = N(L^*) \oplus \overline{R(L)}  $, let $L' = L^*|_{Dom(L^*)\cap \overline{R(L)}}$,     
         $L'\text{ is injective},$ which implies
        $C_1 \text{ is injective}.$ Now, let $g_n \in R(C_1)$ such that $g_n \rightarrow g$ as $n \rightarrow \infty$. 
    Then there exists $h_n \in Dom(L^*)\cap \overline{R(L)}$ such that $g_n = CL^* h_n$. Since $C^{-1}$ is bounded and $\overline{R(L)}$ is closed, there exists $h \in Dom(C)$ such that $CL^* h_n \rightarrow C h$ and $L^* h_n \rightarrow h$. Since $L^*$ is closed and $Dom(L^*) = N(L^*) \oplus Dom(L^*) \cap \overline{R(L)}$, $L'$ is closed, and $R(L') = R(L^*)$. Since, $R(L')$ is closed and $L'$ is injective, there exists $h' \in {Dom(L^*) \cap \overline{R(L)}}$ such that $L^* h_n \rightarrow L^*h'$ and $h_n \rightarrow h'$. Therefore $g = CL^*h'$. Then $g \in R(C_1)$.        
Therefore $\{Lf_n\}_{n \in \mathbb{N}}$ is a lower semi-frame for $\overline{R(L)}$.
\end{proof}
\begin{eg}
      Let $H = \ell_2$ and define $\{f_n\} = \{ne_n\}$, where $\{e_n\}_{n \in \mathbb{N}}$ is the standard orthonormal basis for $H$. Define the operator $$L: \ell_2 \rightarrow \ell_2 \text{ as }L(x_1,x_2,\ldots) = (x_1+x_2,0,x_3,x_4,\ldots).$$ Then, $Dom(L) = H$, and the range of $L$ is $R(L)= span\{e_1,0,e_3,e_4,\ldots\}$. Next the adjoint of $L$ is given by $L^*(x_1,x_2,\ldots) = (x_1,x_1,x_3,x_4,\ldots)$, which implies that the range of $R(L^*) = span\{e_1+e_2,e_3,e_4,\ldots\}$ and $R(L^*) $ is closed,   $Dom(L^*) \cap R(L) = H \cap R(L) = R(L).$ The sequence, $\{Lf_n\} = \{e_1,2e_1,3e_3,\ldots\}$. Now, let $f = (x_1,0,x_3,\ldots) \in R(L)$, $$\sum_{n \in \mathbb{N}} | \langle f,f_n\rangle|^2 = \sum_{n = 3}^{\infty} | nx_n|^2 + 3| x_1|^2 \geq \sum_{n \in \mathbb{N}} | x_n|^2 = \|f\|^2.$$ Therefore, $\{Lf_n\}_{n \in \mathbb{N}}$ is a lower semi-frame for $R(L).$      
      \end{eg}
\begin{prop}
    Let $\{f_n\}_{n \in \mathbb{N}}$ be a sequence in $H$ and $L$ be an operator  with dense domain on $H$ with $R(L^*)=H$. Then $\{Lf_n\}_{n \in \mathbb{N}}$ is a lower semi-frame sequence for $H$ if and only if $\{f_n\}_{n \in \mathbb{N}}$ is a lower semi-frame for $H$.
\end{prop}
\begin{proof}
By assuming that $\{Lf_n\}_{n\in \mathbb{N}}$ is a lower semi-frame sequence, we can conclude that the analysis operator $CL^*|_{Dom(L^*)\cap \overline{R(L)}}$ is injective and $R(CL^*|_{Dom(L^*)\cap \overline{R(L)}})$ $ = R(C)$ is closed. Since $CL^*|_{Dom(L^*)\cap \overline{R(L)}}$ and $L^*|_{Dom(L^*)\cap\overline{R(L)}}$ both are injective, so is $C$. Therefore, $\{f_n\}_{n \in \mathbb{N}}$ is a lower semi-frame for $H.$ \\Conversely, if $\{f_n\}_{n \in \mathbb{N}}$ is a lower semi-frame for $H$ and $R(L^*)=H$, then By the Proposition \ref{e}, $\{Lf_n\}_{n \in \mathbb{N}}$ is a lower semi-frame sequence for $H.$
\end{proof}

      \begin{prop}\label{RFS}
      Let $\{f_n\}_{n \in \mathbb{N}}$ be a Riesz-Fischer sequence for $H$. Let $L: Dom(L)\subseteq H \rightarrow H$ be an operator  with dense domain with $R(L^*) = H$. Then $\{Lf_n\}_{n \in \mathbb{N}}$ is a Riesz-Fischer sequence for $H$.
  \end{prop}
  \begin{proof}
    Given our assumption on $\{f_n\}_{n \in \mathbb{N}}$ we conclude that $C$ is surjective. Since $R(L^*) = H$ it follows that $Dom(CL^*) = Dom(C) \cap R(L^*) = Dom(C)$. Therefore, $R(CL^*) = R(C) = H$. Thus, $\{Lf_n\}$ is a Riesz-Fischer sequence for $H$.    
  \end{proof}
 
\section{Generating sequences using sum}

In this section, we investigate the stability of lower semi-frames and Riesz-Fischer sequences under various operator-induced perturbations in a Hilbert space setting. We present several results that provide sufficient conditions for the preservation of the lower semi-frame property when the original sequence is modified by operators with specific analytical properties. In particular, we show that if $\{f_n\}$ is a lower semi-frame and $L$ is a densely defined operator satisfying certain spectral and range conditions, then the perturbed sequence $\{f_n + Lf_n\}$ remains a lower semi-frame. This is extended to cases where the perturbation involves the sum of an unbounded and a bounded operator, provided that the norm of the bounded operator is controlled relative to the spectral bound of the unbounded operator’s adjoint. Additionally, we establish a similar stability result for Riesz-Fischer sequences and demonstrate that the sum of a lower semi-frame and a Bessel sequence also forms a lower semi-frame under completeness and norm conditions.

\begin{prop}
    Let $\{f_n\}_{n \in \mathbb{N}}$ be a lower semi-frame for $H,$ and let $L$ be an operator  with dense domain on $H$. If $R(I+L)$ and $R(L)$ are dense in $H$, and $\gamma(L^*)> 1 $, then $\{f_n+Lf_n\}_{n \in \mathbb{N}}$ is also a lower semi-frame for $H$.
\end{prop}
\begin{proof}
    The analysis operator for $\{f_n+Lf_n\}_{n \in \mathbb{N}}$ is $C(I+L^*).$ Let $f \in Dom(I+L^*)$, then $
     \| (I+L^*) f \| \geq |  \| L^*f \|-\| f \|  | 
        \geq (\gamma(L^*) - 1) \| f \|$    
    Thus, $ \gamma(I+L^*) = \gamma(L^*) - 1$.
    By the Proposition \ref{closed}, $R(I+L^*)$ is closed and by the Proposition $\ref{a}$, we can conclude that $\{f_n + Lf_n\}$ is a lower semi-frame for $H$.    
    \end{proof}
    
      \begin{cor}
          Let $\{f_n\}_{n \in \mathbb{N}}$ be a lower semi-frame for $H$ and let $L$ be an operator  with dense domain on $H$. If $R(L)$ is dense and $\gamma(L^*)> 1 $, then $\{f_n+Lf_n\}_{n \in \mathbb{N}}$ is a lower semi-frame sequence for $H$.
      \end{cor}
     
          \begin{proof}
 The assumption $\gamma(L^*)>1$ implies that $\| L^*f \| \geq \gamma(L^*) \| f \|$ for all $f \in Dom(L^*)$.
    Let $f \in Dom(I+L^*) \cap N(I+L^*)^\perp$. Then,
        $\| (I+L^*) f \| \geq |  \| L^*f \|-\| f \|  | 
        \geq (\gamma(L^*) - 1) \| f \|$.    
    Therefore, $ \gamma(I+L^*) = \gamma(L^*) - 1$. That is, $R(I+L^*)$ is closed. Our assumption in $L$ shows that $I+L$ is also an operator  with dense domain in $H$. Then the Proposition $\ref{e}$, shows that $\{f_n + Lf_n\}$ is a lower semi-frame sequence for $H$. 
    \end{proof}
    \begin{prop}
        Let $\{f_n\}$ be a lower semi-frame for $H$. Let $L_1 : Dom(L_1) \subseteq H \rightarrow H$ be an operator  with dense domain on $H$ and let $L_2 $ be a bounded operator on $H$, If $R(L_1+L_2)$ is dense ,  and $\gamma(L_1^*) > \| L_2 \| $, then $\{(L_1+L_2)f_n\}_{n \in \mathbb{N}}$ is a lower semi-frame for $H$.
    \end{prop}
    \begin{proof}
       From the hypothesis, it follows that $(L_1+L_2)^*$ is injective.
        Let $f \in Dom((L_1+L_2)^*)$. As $\gamma(L_1^*) > \|L_2\|$ and  $(L_1+L_2)^* \subseteq L_1^*+L_2^*$, we have $\|(L_1+L_2)^*f\| = \|L_1^*f+L_2^*f\| \geq | \|L_1^*f\|-\|L_2^*f\|| \geq (\gamma(L_1^*) - \|L_2\|) \|f\|$. Therefore, $R(L_1+L_2)^*$ is closed. this means that $\{(L_1+L_2)f_n\}_{n \in \mathbb{N}}$ is a lower semi-frame for $H$.      
    \end{proof}
  
   \begin{cor}
          Let $\{f_n\}_{n \in \mathbb{N}}$ be a Ries-Fischer sequence for $H$. Let $L$ be an operator  with dense domain on $H$. If $R(L^*) = H$  and $\gamma(L^*)> 1 $, then $\{f_n+Lf_n\}_{n \in \mathbb{N}}$ is a Riesz-Fischer sequence for $H$.
      \end{cor}
     
          \begin{proof}
 From the hypothesis, we observe that $L$ is injective and $(I+L)$ is also injective. Consequently, $R(I+L^*)$ is dense in $H$. Consider the sequence $\{f_n+Lf_n\}_{n \in \mathbb{N}}$ in $H.$  
  Let $f \in Dom(I+L^*)$. Then, $
     \| (I+L^*) f \| \geq |  \| L^*f \|-\| f \|  | 
        \geq (\gamma(L^*) - 1) \| f \|.$     
    Therefore, $ \gamma(I+L^*) = \gamma(L^*) - 1$.
    By the Proposition \ref{closed}, $R(I+L^*)$ is closed. Since $R(I+L^*)$ is dense, $R(I+L^*) = H$. By the Proposition \ref{RFS}, $\{f_n+Lf_n\}_{n\in \mathbb{N}}$ is a Riesz-Fischer sequence for $H$. 
    \end{proof}

    \begin{prop}
        Let  $\{f_n\}_{n \in \mathbb{N}}$ be a lower semi-frame for $H$ with analysis operator $C_1$ and  lower frame bound  $\alpha,$ and let $\{g_n\}_{n \in \mathbb{N}}$ be a Bessel sequence in $H$ with analysis operator $C_2$ and Bessel bound $\beta$. If $\{f_n+g_n\}_{n \in \mathbb{N}}$ is complete in $H$ and  $\sqrt{\alpha} > \sqrt{\beta}$, then $\{f_n+g_n\}_{n \in \mathbb{N}}$ is a lower semi-frame for $H$.
        \end{prop}
        \begin{proof}
           Let the sequence $\{f_n + g_n\}_{n \in \mathbb{N}}$ in $H$. The analysis operator for this sequence is given by  $\{\langle f, f_n+g_n \rangle\}_{n\in \mathbb{N}} = \{\langle f, f_n \rangle\}_{n\in \mathbb{N}} + \{\langle f, g_n \rangle\}_{n\in \mathbb{N}} = C_1+C_2$.
            Since $\{f_n+g_n\}_{n \in \mathbb{N}}$ is complete, it follows that $C_1+C_2$ is injective. If $\sqrt{\alpha} > \sqrt{\beta}$, for any $f\in H,$  $ \|(C_1+C_2)f\| \geq | \|C_1f\|-\|C_2f\|| \geq \|C_1f\|-\|C_2f\| \geq (\sqrt{\alpha} -\sqrt{\beta})\|f\|$. Therefore, $\{f_n+g_n\}_{n \in \mathbb{N}}$ is a lower semi-frame for $H$.              
        \end{proof}
    
    \section{Generating sequences using direct sum}
    We explored the direct sum of sequences as a generalization of the set theoretic sum. We attain necessary and sufficient conditions for preserving the frame-like structure of the sequences under direct sum. Through these results we get that if $H$ can be written as the direct sum of subspaces of $H$ and the sequences in the subspaces satisfy certain conditions, then the frame-like structure of sequences can be hereditary.  We begin this section by considering sequences $\{f_n\}_{n \in \mathbb{N}} \subset H_1 $ and $\{g_n\}_{n \in \mathbb{N}} \subset H_2$ with analysis operators $C_1$ and $C_2$, respectively. Then the analysis operator for $\{f_n \oplus g_n\}_{n \in \mathbb{N}} \subset H_1 \oplus H_2$ is $C: Dom(C) \subseteq H_1 \oplus H_2 \rightarrow \ell_2$ defined by $$ C(f\oplus g) = \{\langle f\oplus g , f_n\oplus g_n \rangle\}_{n \in \mathbb{N}}=\{\langle f, f_n \rangle_{H_1} + \langle g, g_n \rangle_{H_2} \}_{n \in \mathbb{N}} = C_1f+C_2g.$$
 Then, $Dom(C) = Dom(C_1) \oplus Dom(C_2)$ and $R(C) = R(C_1)+R(C_2).$ \\ \hspace*{0.5cm}We are adapting the definitions \ref{def} of certain frames characterized by Deguang and Larson in \cite{Dsum} for general sequences and defining in the following manner.
\begin{defn}
    Let $\{f_n\}_{n \in \mathbb{N}} $ and $\{g_n\}_{n \in \mathbb{N}}$ be sequences in $H_1$ and $H_2$ respectively with analysis operators $C_1$ and $C_2$ respectively. Then
    \begin{enumerate}[(i)]
        \item $\{f_n\}_{n \in \mathbb{N}} $ and $\{g_n\}_{n \in \mathbb{N}}$ are said to be strongly complementary if $R(C_1)\oplus ^\perp R(C_2) = \ell_2$.
        \item $\{f_n\}_{n \in \mathbb{N}} $ and $\{g_n\}_{n \in \mathbb{N}}$ are said to be complement if $R(C_1) \cap R(C_2) = \{0\}$ and $R(C_1)+R(C_2) = \ell_2$.
        \item $\{f_n\}_{n \in \mathbb{N}} $ and $\{g_n\}_{n \in \mathbb{N}}$ are said to be strongly disjoint if $R(C_1) \perp R(C_2)$.
        \item $\{f_n\}_{n \in \mathbb{N}} $ and $\{g_n\}_{n \in \mathbb{N}}$ are said to be disjoint if $R(C_1) \cap R(C_2) = \{0\}$.
    \end{enumerate}
\end{defn} 

\begin{prop}\label{4.2}
    Let $\{f_n\}_{n \in \mathbb{N}} \subset H_1$ and $\{g_n\}_{n\in \mathbb{N}} \subset H_2$ be disjoint sequences. Then $\{f_n \oplus g_n\}_{n \in \mathbb{N}}$ is a complete sequence for $H_1 \oplus H_2$ if and only if $\{f_n\}_{n \in \mathbb{N}}$ and $\{g_n\}_{n \in \mathbb{N}}$ are complete sequences for $H_1$ and $H_2$ respectively.
    \end{prop}
    \begin{proof}
        The completeness of the sequence $\{f_n \oplus g_n\}_{n \in \mathbb{N}}$ implies that the analysis operator $C$ is injective, that is $N(C) = \{0\}$. Let $f \in N(C_1)$ and $g \in N(C_2)$,
        \begin{center}
            $\{\langle f \oplus g, f_n \oplus g_n \rangle \}_{n \in \mathbb{N}} = \{\langle f, f_n \rangle\}_{n \in \mathbb{N}} + \{\langle g, g_n \rangle\}_{n \in \mathbb{N}} = 0$.
        \end{center}
        This implies that $f \oplus g \in N(C)$ and $f \oplus g = 0$. Therefore, $C_1$ and $C_2$ are injective, which in turn gives $\{f_n\}_{n \in \mathbb{N}}$ and $\{g_n\}_{n \in \mathbb{N}}$ are complete in $H_1$ and $H_2$ respectively.\\
        \hspace*{0.5cm}Conversely, assume that $f \oplus g \in N(C).$ Then,
        $\{\langle f \oplus g, f_n \oplus g_n \rangle \}_{n \in \mathbb{N}} = \{\langle f, f_n \rangle\}_{n \in \mathbb{N}}$ $ + \{\langle g, g_n \rangle\}_{n \in \mathbb{N}} = 0.$ Since the sequences $\{f_n\}_{n \in \mathbb{N}}$ and $\{g_n\}_{n\in \mathbb{N}}$ are complete in their respective spaces, it follows that $\{\langle f, f_n \rangle_{H_1}\}_{n \in \mathbb{N}} \in R(C_1) \cap R(C_2) = \{0\}$. Thus $f \in N(C_1) = \{0\}$, and similarly, $g = 0$. This confirms that $C$ is injective  and therefore, the sequence $\{f_n \oplus g_n\}_{n \in \mathbb{N}}$ is complete in $H_1 \oplus H_2$.       
    \end{proof}
     \begin{cor}
        Let $K_1$ and $K_2$ be closed subspaces of $H$ such that $K_1 \oplus ^\perp  K_2 = H$. If $\{f_n\}_{n \in \mathbb{N}}$ and $\{g_n\}_{n \in \mathbb{N}}$ are complete sequences for $K_1$ and $K_2$ respectively such that $\{f_n\}_{n \in \mathbb{N}}$ and $\{g_n\}_{n \in \mathbb{N}}$ are disjoint, then $\{f_n \oplus g_n\}_{n \in \mathbb{N}}$ is a complete sequence for $H$.
    \end{cor}
    \begin{proof}
        Since $K_1$ and $K_2$ are closed subspaces of $H$, $K_1$ and $K_2$ are complete. Consider $H_1 = K_1$ and $H_2 = K_2$ in Proposition \ref{4.2}, then $H_1 \oplus H_2 = H$. Therefore, $\{f_n \oplus g_n\}_{n \in \mathbb{N}}$ is a complete sequence for $H.$        
    \end{proof}
    \begin{prop}\label{4.3}
        Let $\{f_n\}_{n \in \mathbb{N}} \subset H_1$ and $\{g_n\}_{n \in \mathbb{N}} \subset H_2$ be lower semi-frames. If $\{f_n\}_{n \in \mathbb{N}}$ and $\{g_n\}_{n \in \mathbb{N}}$ are disjoint sequence, then $\{f_n \oplus g_n\}_{n \in \mathbb{N}} \subset H_1\oplus H_2$ is a lower semi-frame.
    \end{prop}
    \begin{proof}
       Let $C_1$ and $C_2$ be the analysis operators for $\{f_n\}_{n \in \mathbb{N}}.$ and $\{g_n\}_{n \in \mathbb{N}}$ respectively.  By the Proposition \ref{4.2}, $\{f_n \oplus g_n\}_{n \in \mathbb{N}}$ is a complete sequence in $ H_1\oplus H_2$. Now, let $\{x_n \oplus y_n\}_{n \in \mathbb{N}} \in R(C_1)+R(C_2)$ such that $\{x_n \oplus y_n\}_{n \in \mathbb{N}} \rightarrow x \in \ell_2$. Since $R(C_1)$ is closed, $\ell_2 = R(C_1) \oplus ^\perp  R(C_1)^{\perp}$. Then there exists $x_1' \in R(C_1)$ and $x_2' \in R(C_1)^{\perp}$ such that $x = x_1'+x_2'$.\\ \hspace*{0.5cm} Since $\{x_n \oplus y_n\}_{n \in \mathbb{N}} \rightarrow x= x_1'+x_2'$ as $n \rightarrow \infty$, $x_n \rightarrow x_1'$ and $y_n \rightarrow x_2'$. Since $R(C_2)$ is closed, $\ell_2 = R(C_2) \oplus ^\perp  R(C_2)^{\perp}$. Then there exists $x_1 \in R(C_2)$ and $x_2 \in R(C_2)^{\perp}$ such that $x = x_1+x_2$. Since $\{x_n \oplus y_n\}_{n \in \mathbb{N}} \rightarrow x= x_1+x_2$ as $n \rightarrow \infty$, $x_n \rightarrow x_2$ and $y_n \rightarrow x_1$. By algebra of limits, $x_n \rightarrow x_1'=x_2$ and $y_n \rightarrow x_2'=x_1$ as $n \rightarrow \infty$. Therefore, $\{x_n \oplus y_n\}_{n \in \mathbb{N}} \rightarrow x= x_1'+x_2' \in R(C_1)+R(C_2)$ as $n \rightarrow \infty$. From this, we can conclude that $\{f_n \oplus g_n\}_{n \in \mathbb{N}}$ is a lower semi-frame for $H_1 \oplus H_2$.
    \end{proof}
    \begin{prop}\label{4.5}
       Let $\{f_n\}_{n \in \mathbb{N}} \subset H_1$ and $\{g_n\}_{n\in \mathbb{N}} \subset H_2$ be strongly disjoint sequences. Then $\{f_n \oplus g_n\}_{n \in \mathbb{N}}$ is lower semi-frame if and only if $\{f_n\}_{n \in \mathbb{N}}$ and $\{g_n\}_{n \in \mathbb{N}}$ are lower semi-frames for $H_1$ and $H_2$ respectively. 
    \end{prop}
    \begin{proof}
     Since the sequence $\{f_n \oplus g_n\}_{n \in \mathbb{N}}$ forms a lower semi frame for $H_1 \oplus H_2,$ the range $R(C) = R(C_1)+R(C_2)$ is closed. By Proposition \ref{4.2}, it follows that $\{f_n\}_{n \in \mathbb{N}}$ and $\{g_n\}_{n \in \mathbb{N}}$ are complete sequences in $H_1$ and $H_2$ respectively. Now, let $\{x_n\}_{n \in \mathbb{N}} \in R(C_1)$ be a sequence such that $\{x_n\} \rightarrow x$ as $n \rightarrow \infty.$ We can express each term as $x_n $ as $x_n +0$. Then $\{x_n\}_{n \in \mathbb{N}} \in R(C_1)+R(C_2)$ and $\{x_n\} \rightarrow x$ as $n \rightarrow \infty$, which implies that $x \in R(C_1)+R(C_2).$ That is, there exists $x_1' \in R(C_1)$ and $x_2' \in R(C_2)$ such that $x = x_1'+x_2'$. Therefore $\{x_n\} \rightarrow x_1'+x_2'$ as $n \rightarrow \infty$.\\ \hspace*{0.5cm} Since $R(C_1) $ and $R(C_2)$ are orthogonal, we have $\{x_n\} \rightarrow x_1'$ as $n \rightarrow \infty$. This implies that $x_2' = 0$ and $x = x_1' \in R(C_1).$ Therefore, $R(C_1)$ is closed. Similarly, we can show that $R(C_2)$ is  closed. Thus, $\{f_n\}_{n \in \mathbb{N}}$ and $\{g_n\}_{n \in \mathbb{N}}$ are lower semi-frames for $H_1$ and $H_2$ respectively. \\ \hspace*{0.5cm}Conversely, from the Proposition \ref{4.3}, we conclude that $\{f_n \oplus g_n\}_{n \in \mathbb{N}}$ is a lower semi-frame for $H_1 \oplus H_2$.       
    \end{proof}
    \begin{cor}
        Let $K_1$ and $K_2$ be closed subspaces of $H$ such that $K_1 \oplus ^\perp  K_2 = H$. If $\{f_n\}_{n \in \mathbb{N}}$ and $\{g_n\}_{n \in \mathbb{N}}$ are lower-semi frames for $K_1$ and $K_2$ respectively such that $\{f_n\}_{n \in \mathbb{N}}$ and $\{g_n\}_{n \in \mathbb{N}}$ are strongly disjoint, then $\{f_n \oplus g_n\}_{n \in \mathbb{N}}$ is a lower semi-frame for $H$.
    \end{cor}
    \begin{proof}
        Since $K_1$ and $K_2$ are closed subspaces of $H$, $K_1$ and $K_2$ are complete. Consider $H_1 = K_1$ and $H_2 = K_2$ in Proposition \ref{4.5}, then $H_1 \oplus H_2 = H$. Therefore, $\{f_n \oplus g_n\}_{n \in \mathbb{N}}$ is a lower semi-frame for $H.$        
    \end{proof}
     \begin{eg}
        Let $\{e_n\}_{n \in \mathbb{N}}$ be the standard orthonormal basis for $\ell_2.$ Let $H_1 = span\{e_{2n} : n \in \mathbb{N}\}$ and $H_2 = span\{e_{2n+1} : n \in \mathbb{N}\}$. Consider $\{f_n\}_{n \in \mathbb{N}} = \{2ne_{2n}\}_{n \in \mathbb{N}}$ and $\{g_n\}_{n \in \mathbb{N}} = \{(2n+1)e_{2n+1}\}_{n \in \mathbb{N}}$. The analysis operator for $\{f_n\}_{n \in \mathbb{N}}$ is $C_1(x_1,x_2,\ldots) = (x_1,0,3x_3,\ldots)$ and $R(C_1) = H_1$ and the analysis operator for $\{g_n\}_{n \in \mathbb{N}}$ is $C_2(x_1,x_2,\ldots) = (0,2x_2,0,4x_4,\ldots)$ and $R(C_2) = H_2$. Therefore $R(C_1) \perp R(C_2) .$ Since $R(C_1)$ and $R(C_2)$ are closed, $\{f_n \}_{n \in \mathbb{N}}$ and $\{g_n\}_{n \in \mathbb{N}}$ are lower semi-frames for $H_1$ and $H_2$ respectively. Then, $\{f_n \oplus g_n\}_{n \in \mathbb{N}} = \{ne_{n}\}_{n \in \mathbb{N}}$ which is a lower semi-frame for $\ell_2$. 
    \end{eg}
    
    \begin{prop}
        Let $\{f_n\}_{n \in \mathbb{N}} \subset H_1$ be a Riesz-Fischer sequence and $\{g_n\}_{n \in \mathbb{N}} \subset H_2$ be any sequence. Then $\{f_n \oplus g_n\}_{n \in \mathbb{N}}$ is a Riesz-Fischer sequence for $H_1 \oplus H_2$ 
    \end{prop}
    \begin{proof}
       Let $C_1$ and $C_2$ be the analysis operators for $\{f_n\}_{n \in \mathbb{N}}.$ and $\{g_n\}_{n \in \mathbb{N}}$ respectively. From our assumption, $R(C_1) = \ell_2.$ Assume that the analysis operator for $\{f_n \oplus g_n\}_{n \in \mathbb{N}}$ is $C$. Then $R(C) = R(C_1)+R(C_2) = \ell_2+R(C_2) =\ell_2.$ This concludes that $\{f_n \oplus g_n\}_{n \in \mathbb{N}}$ is a Riesz-Fischer sequence for $H_1 \oplus H_2.$       
    \end{proof}

\begin{center}
	\textbf{Conclusion}
\end{center}
 In this paper, we established conditions on a class of unbounded operators that preserve the frame-like properties of sequences in Hilbert spaces. The primary motivation for this study stemmed from the question of whether a sequence in a Hilbert space can inherit the same structural properties as sequences in its subspaces. This inquiry led to an exploration of the direct sum of sequences, where we derived necessary and sufficient conditions for maintaining frame-like structures.
\begin{center}
	\textbf{Acknowledgements}
\end{center}

\noindent The first author gratefully acknowledges the National Institute of Technology Karnataka (NITK), Surathkal, for their financial support. The second and third authors thank ANRF (SERB), DST, Government of India (TAR/2022/000219) for the TARE Fellowship Grant.
\addcontentsline{toc}{section}{References}

\end{document}